\documentclass[12pt,oneside,reqno]{amsart}
\usepackage{amsmath,amsfonts,amssymb,stmaryrd}

\usepackage{color}

\usepackage{tikz}
\usetikzlibrary{calc}
\usepackage{subfig}

\newtheorem{theorem}{Theorem}[section]
\newtheorem{corollary}[theorem]{Corollary}
\newtheorem{lemma}[theorem]{Lemma}

\theoremstyle{definition}

\theoremstyle{remark}
\newtheorem{remark}[theorem]{Remark}

\numberwithin{equation}{section}
\def\u{{\bf u}}
\def\O{\Omega}

\def\R{\mathbb{R}}

\def\d{\mbox{div\,}}
\def\sop{\mbox{supp\,}}
\def\G{{\bf G}}
\def\ve{\varepsilon}
\def\vphi{\varphi}

\def\n{\nabla}
\def\a{\alpha}

\begin{document}

	\date{}

\title{Solutions of the divergence equation in Hardy and lipschitz spaces}

\author{Mar\'\i a E. Cejas}
\address{IMAS (UBA-CONICET) 
	Ciudad Universitaria\\
	(1428) Ciudad Aut\'onoma de Buenos Aires\\  
	Argentina;
	e-mail: ecejas@mate.unlp.edu.ar}

\author{Ricardo G. Dur\'an}
\address{IMAS (UBA-CONICET) and Departamento de Matem\'atica\\
Facultad de Ciencias Exactas y Naturales\\
Universidad de Buenos Aires\\
Ciudad Universitaria\\
(1428) Ciudad Aut\'onoma de Buenos Aires\\
Argentina;	e-mail: rduran@dm.uba.ar.}



\begin{abstract}
Given a bounded domain $\O$ and $f$ of zero integral, the existence of a vector fields $\u$ vanishing on $\partial\O$ and satisfying $\d\u=f$ has been widely studied because of its connection with many important problems. It is known that for $f\in L^p(\O)$, $1<p<\infty$, there exists a solution $\u\in W^{1,p}_0(\O)$, and also that an analogous result is not true for $p=1$ or $p=\infty$. The goal of this paper is to prove results for Hardy spaces when $\frac{n}{n+1}<p\le 1$, and in the other limiting case, for bounded mean oscillation and Lipschitz spaces. As a byproduct of our analysis we obtain a Korn inequality for vector fields in Hardy-Sobolev spaces.
\end{abstract}

\maketitle
\section{Introduction}
For a bounded domain $\O\subset\R^n$, $n\ge 2$, and  $1\le p\le\infty$, we will use the standard notation $W^{1,p}(\O)$, for the Banach space of functions in $L^p(\O)$ with first derivatives also in $L^p(\O)$, and $W^{1,p}_{0} (\O)$ for the subspace of functions vanishing on the boundary of $\O$.
Also we denote by $L^p_0(\O)$ the subset of $L^p(\O)$ of functions such that $\int_{\O}f=0$.

Given $f\in L^p_0(\O)$, the existence of solutions of
\begin{equation}
\label{div}
\left\{
\begin{array}{rl}
\d\u = f &\ \ \mbox{in} \ \O\\
\u=0& \ \ \mbox{on} \ \partial\O
\end{array}
\right.
\end{equation}
is a classic problem that has been widely studied because of it connections with many other results which are fundamental in the analysis of partial differential equations arising in classic mechanics. 

For example, when $p=2$ the existence of $\u\in W_0^{1,2}(\O)^n$ satisfying \eqref{div} as well as $\|\u\|_{W^{1,2}(\O)} \le C\|f\|_{L^2(\O)}$, is a key result to prove the well posedness of the classic Stokes equations. Consequently, several arguments have been introduced to prove this result under different assumptions on the domain $\O$ (see for example \cite{AD} for a historical account and references). 

The extension of this result to the case $1<p<\infty$ was
proved by Bogovski\v\i{} \cite{B} by means of an explicit integral operator giving a solution of \eqref{div} for domains which are star-shaped with respect to a ball. Then, he extended the result to bounded Lipschitz domains using that they can be written as the union of a finite number of star-shaped domains.
Precisely, Bogovski\v\i {} proved the following:

\begin{theorem}
\label{teorema en Lp}
Let $\O \subset \R^n$, $n \ge 2$, be a bounded Lipschitz  domain. Given $1<p<\infty$ and $f \in L^p_0(\O)$, there exists $\u\in W^{1,p}_0(\O)^n$ such that 
\begin{equation*}
\d\u  =  f \ \mbox{in} \ \O
\end{equation*}
and
\begin{equation}
\label{cota en Lp}
\|\u\|_{W^{1,p}(\O)}\le C\|f\|_{L^p(\O)}
\end{equation}
where $C= C(\O, p,n)$. 
\end{theorem}
After Bogovski\v\i's paper his operator has been widely analyzed
by many authors  and several generalizations have been obtained (see for example \cite{ADM,CM,DM,GHH,S}). 

It is known that Theorem \ref{teorema en Lp} cannot be
generalized to the limiting cases
$p=1$ and $p=\infty$. Actually, the result is false even 
without imposing boundary conditions. Indeed
$\d\u=f$ has, in general, no 
solution in $W^{1,1}(\O)$ (resp. $W^{1,\infty}(\O)$) for 
$f\in L^1(\O)$ (resp. $f\in L^\infty(\O)\cap C^0(\O))$. Many authors have considered these problems (we refer to \cite{DFT} and references therein).

To illustrate we give here a simple argument showing that
the result is not valid for $p=1$.
Consider for example $n=2$ and let $D$ be the unit disk centered at the origin. If $\phi$ is a harmonic function such that $\int_D\phi=0$ in $D$ we denote by $\psi$ its harmonic conjugate, namely, 
$\phi_x=\psi_y$ and $\phi_y=-\psi_x$ in $D$.
Suppose that for any $f\in L^1(D)$ there exists $\u\in W^{1,1}_0(D)^2$ satisfying
$\d\u=f-\bar f$ (where $\bar f$ is the average of $f$) and $\|\u\|_{W^{1,1}(D)}\le C\|f\|_{L^1(D)}$, then
$$
\begin{aligned}
\int_D \phi \, f&=\int_D \phi \, \d\u=-\int_D (\phi_x u_1+\phi_y u_2)
=-\int_D (\psi_y u_1-\psi_x u_2)
\\
&=\int_D (\psi u_{1y}-\psi u_{2x})\le C \|\psi\|_{L^\infty}\|\u\|_{W^{1,1}}
\le C \|\psi\|_{L^\infty}\|f\|_{L^1}
\end{aligned}
$$ 
and therefore, we conclude by duality that 
$\|\phi\|_{L^\infty(D)}\le \|\psi\|_{L^\infty(D)}$, which is not true if we take, for example, $\phi$ as the real part
of $\log(1+z)$.

Consequently, a natural question is whether the result of Theorem \ref{teorema en Lp} can be extended to some values of $p\le 1$ replacing the classic Sobolev spaces by Hardy-Sobolev ones. 

For the first question a positive answer was given in \cite{ChJY} for $\frac{n}{n+1}<p\le 1$. Our first goal is to give a simpler proof of the results of that paper. Indeed, we will show that Bogovski\v\i's operator is bounded in Hardy-Sobolev spaces. We also generalize the results to weighted spaces.

Analogously, at the other limiting case, one can ask whether similar results hold for $BMO$ and Lipschitz-$\a$ spaces. To answer this question is our second goal. 

Finally, we prove a Korn inequality in Hardy-Sobolev spaces for $\frac{n}{n+1}<p\le 1$. For $p=1$, the so called first case of Korn inequality is proved in \cite{SS}, where the authors leave as an open question whether the second case is valid. We give a positive answer to this question proving the general form of Korn inequality from which both first and second cases can be derived by usual compactness arguments.

\section{Hardy spaces and their duals}
\label{espacios de hardy}

In this section we recall the definitions of the functional spaces that we will work with.

To recall the definition of the Hardy spaces 
we first introduce the Fefferman-Stein maximal function
and its local version due to Goldberg (see \cite{G,St}).

For a fixed $\vphi\in\mathcal S$ such that $\int\vphi=1$, we use the standard notation
$\vphi_t(x)=t^{-n}\vphi(x/t)$, and define, for $f\in\mathcal S'$,  
$$
M_{\vphi}f(x)=\sup_{t>0} |(f*\vphi_t)(x)|
$$
and
$$
m_{\vphi}f(x)=\sup_{0<t<1} |(f*\vphi_t)(x)|.
$$
Then, for $0<p<\infty$, the global Hardy spaces are given by
$$
H^p(\R^n)=\{f\in\mathcal S'\colon M_{\varphi}f\in L^{p}(\R^n)\},
$$
while the local ones, by 
$$
h^p(\R^n)=\{f\in\mathcal S'\colon m_{\varphi}f\in L^{p}(\R^n)\}.
$$
These are complete spaces with the norms (or quasi-norms if $0<p<1$) defined as
$\|f\|_{H^p(\R^n)}=\|M_{\vphi}f\|_{L^p(\R^n)}$ and 
$\|f\|_{h^p(\R^n)}=\|m_{\vphi}f\|_{L^p(\R^n)}$ respectively.

Clearly, $H^p(\R^n)\subset h^p(\R^n)$. The  $h^p(\R^n)$ spaces are more convenient in many applications because they have the important property:

\begin{equation}
\label{eta por f esta en hp}
\eta\in C_0^\infty(\R^n), f\in h^p(\R^n) 
\implies \eta f\in h^p(\R^n).
\end{equation}

To define the Hardy spaces in a domain
$\O\subset\R^n$ two natural possibilities have been considered, namely,
$$
h_z^p(\O)=\{f\in h^p(\R^n)\colon \sop f \subseteq\overline\O\}
$$
and the restricted space,
$$
h_r^p(\O)=\{f=F|_\O\colon F\in h^p(\R^n)\}
$$
with the norm (or quasi-norm if $p<1$) given by 
$$
\|f\|_{h_r^p(\O)}=\inf\{\|F\|_{h^p(\R^n)}
\colon f=F|_\O\}.
$$
We will work also with the Hardy-Sobolev spaces 
$$
h^{1,p}_{z,0}(\O)=\{f \in h_z^p(\O): \nabla f \in h^p_z(\O)^n\}
$$
and
$$
h^{1,p}_r(\O)=\{f \in h_r^p(\O): \nabla f \in h^p_r(\O)^n\}
$$ 
with the norms (or quasi-norms when $0<p<1$) defined as
$\|f\|_{h^{1,p}_z(\O)}=\|f\|_{h^p_z(\O)}+\|\nabla f\|_{h^p_z(\O)}$
and
$\|f\|_{h^{1,p}_r(\O)}=\|f\|_{h^p_r(\O)}+\|\nabla f\|_{h^p_r(\O)}$ respectively.

Now we introduce the homogeneous Lipschitz spaces, for $0<\a<1$,
$$
\dot\Lambda_\a(\R^n)
=\left\{ f\in C(\R^n)\colon \|f\|_{\dot\Lambda_\a}
:= \sup_{x\not= y} \frac{|f(x)-f(y)|}{|x-y|^\a}<\infty
\right\}
$$
and the inhomogeneous ones,
$$
\Lambda_\a(\R^n)
=\left\{ f\in C(\R^n)\colon \|f\|_{\Lambda_\a}
:= \|f\|_{\infty}+\|f\|_{\dot\Lambda_\a}<\infty
\right\}.
$$
Finally, we recall the bounded mean oscillation spaces. Define
$$
\|f\|_{BMO}
:= \sup_{x\in\R^n, r>0} \frac{1}{|B(x,r)|}\int_{B(x,r)} |f(y)-f_B|\,dy
$$
and
$$
\|f\|_{bmo}
:= \sup_{x\in\R^n, 0<r\le 1} \frac{1}{|B(x,r)|}\int_{B(x,r)} |f(y)-f_B|\,dy
+\sup_{x\in\R^n, r>1} \frac{1}{|B(x,r)|}\int_{B(x,r)} |f(y)|\,dy
$$
where $f_B$ denotes the average of $f$ over $B(x,r)$,
and the corresponding spaces
$$
BMO(\R^n)
=\left\{f\in L^1_{loc}(\R^n)\colon \|f\|_{BMO}<\infty\right\}
$$
and
$$
bmo(\R^n)
=\left\{ f\in L^1_{loc}(\R^n)\colon \|f\|_{bmo}<\infty\right\}.
$$
It is known \cite{G} that, for $\frac{n}{n+1}<p<1$,
$$
(H^p(\R^n))^*=\dot\Lambda_{n(1/p-1)}(\R^n)
\ \ \mbox{and} \ \
(h^p(\R^n))^*=\Lambda_{n(1/p-1)}(\R^n)
$$ 
while
$$
(H^1(\R^n))^*=BMO(\R^n)
\ \ \mbox{and} \ \
(h^1(\R^n))^*=bmo(\R^n).
$$

\section{Solvability of the divergence in Hardy spaces}
\label{divergencia en hardy}
	
We want to show that, for a bounded Lipschitz domain $\O$ and $p\in(\frac{n}{n+1},1]$, given $f\in h^p_z(\O)$ there exists a solution of $\d\u=f$ such that $\u\in h^{1,p}_{z,0}(\O)^n$.
 
First we consider star-shaped domains and afterwards extend the result for Lipschitz domains using known arguments.

Let $\O$ be a bounded domain which is star-shaped with respect to a ball $B\subset\O$, and choose $\omega\in C_0^\infty(B)$ such that $\int_{B}\omega=1$. Given $f\in L^1(\O)$ we define
\begin{equation}
\label{def de u}
\u(x)=\int_\O \G (x,y) f(y) dy 
\end{equation}
where
\begin{equation}
\label{definicion de G}
\G(x,y)= \int_0^1\frac{(x-y)}{s}
\,\omega\left(y+\frac{x-y}{s}\right)\frac{ds}{s^n}.
\end{equation}

Then, it was proved by Bogovski\v\i{} \cite{B}, that $\u$ is a solution of problem \eqref{div} whenever $\int_\O f=0$, and
that it satisfies \eqref{cota en Lp} 
(see details, for example, in \cite{AD,G}). 

Taking derivatives in \eqref{def de u}, a standard argument
 (details can be seen in \cite{Ga,AD}) gives, for $x\in\O$,
\begin{equation}
\label{derivadas de u}
\frac{\partial u_i}{\partial x_j}= T_{ij}f+\omega_{ij}f
\end{equation}
where
\begin{equation}
\label{omegaij}
\omega_{ij}(x)=\int \frac{z_iz_j}{|z|^2}\omega(x+z)\,dz
\end{equation}
and
\begin{equation}
\label{parte singular}
T_{ij}f(x)=\lim_{\ve\rightarrow 0}\int_{|x-y|>\ve}
\frac{\partial G_i}{\partial x_j}(x,y)f(y)\,dy.
\end{equation}
In \eqref{omegaij}, as in the rest of the paper, we omit
the domain of integration when it is the whole space $\R^n$.

Since $\omega_{ij}$ is a smooth function, the more complicated part is to show the continuity of the singular integral operator $T_{ij}$.

Up to now we have considered functions supported in $\overline\O$. In order to apply the theory of singular integrals in $\R^n$ we can multiply the kernel by  regularized characteristic functions.
Let $\widetilde B$ be
a ball such that $\overline\O\subset\widetilde B$ and take $\chi\in C_0^\infty(\widetilde B)$ such that $\chi\equiv 1$ in $\O$. 
Then, we modify the definition \eqref{parte singular} and  consider, for smooth functions $f$ defined in $\R^n$,
the operator (that we keep calling $T_{ij}$ to simplify notation) defined by,

\begin{equation}
\label{operador general}
T_{ij}f(x)=\lim_{\ve \rightarrow 0}\int_{|x-y|>\ve}
\chi(x) \chi(y) \frac{\partial G_i}{\partial x_j}(x,y) f(y)\, dy.
\end{equation}
Clearly \eqref{derivadas de u} is valid in $\O$ with this new operator.

Introducing now $\eta_i(y,z):=z_i\omega(y+z)$ we have that
$$
\frac{\partial G_i}{\partial x_j}(x,y)
=\int_0^1\frac{\partial \eta_i}{\partial z_j}\left(y,\frac{x-y}{s}\right)\frac{ds}{s^{n+1}}
$$
and so,
$$
T_{ij}f(x)=\lim_{\ve \rightarrow 0}\int_{|x-y|>\ve}
\chi(x) K(y,x-y) f(y)\, dy
$$
where
\begin{equation}
\label{def de K}
K(y,z)=
\int_0^1 \chi(y)\frac{\partial \eta_i}{\partial z_j}\left(y,\frac{z}{s}\right)\frac{ds}{s^{n+1}}.
\end{equation}

To prove the continuity of this operator
we will make use of the following known result.
We will use the standard notation $A\lesssim B$ which means that the quantity $A$ is bounded by a constant times $B$.

\begin{theorem} 
\label{teo de Komori Hp}
Let 
$$
Tf(x)=\lim_{\ve \rightarrow 0}\int_{|x-y|>\ve} N(x,y)f(y)\,dy
$$
be a singular integral operator which is bounded in $L^2$. 
Assume that $T$ satisfies the following conditions:
	
\begin{equation}
\label{K11}
|N(x,y)|\lesssim \min\left\{\frac{1}{|x-y|^n},\frac{1}{|x-y|^{n+1}}\right\}    \end{equation}
	
\begin{equation}
\label{K12}
2|y-\bar y| < |x - y| \implies |N(x, y) - N(x, \bar y)| \lesssim \frac{|y-\bar y|}{|x-y|^{n+1}}.
\end{equation}
Then, for $0<\a<1$,
\begin{equation}
\label{teo de dhz}
T^*1\in\dot\Lambda_\a(\R^n) 
\implies T \, \mbox{is bounded on} \ \  h^p(\R^n) \ \mbox{for all} \ \ \frac{n}{n+\a}<p\le 1.
\end{equation}
 \end{theorem} 

\begin{proof} Replacing hypothesis \eqref{K11} by the weaker one 
$|N(x,y)|\lesssim |x-y|^{-n}$,
Komori proved in \cite[Theorem 2]{K1} that $T$ is bounded
from $H^p(\R^n)$ to $h^p(\R^n)$. His arguments were modified
in \cite{DHZ} to show that, under the stronger decay of $N(x,y)$ assumed in \eqref{K11}, \eqref{teo de dhz} holds.
\end{proof}

We will call $d$ the diameter of $\widetilde B$.

\begin{lemma}
\label{cotas del nucleo}
The kernel $N(x,y):=\chi(x) K(y,x-y)$ satisfies \eqref{K11} and \eqref{K12}.
\end{lemma}

\begin{proof}
It follows from its definition 
that $K(y,z)$ vanishes for $y$ outside 
$\widetilde B$. Then, recalling that $\eta_i(y,z):=z_i\omega(y+z)$ and that
$\sop\omega\subset\O$, it follows that the integrand in \eqref{def de K} vanishes unless $|z|\le sd$.
Therefore, we can rewrite the kernel as,
$$
K(y,z)=
\int_{|z|/d}^1\chi(y)\frac{\partial \eta_i}{\partial z_j}\left(y,\frac{z}{s}\right)\frac{ds}{s^{n+1}}
$$
and so
$$
|K(y,x-y)|\lesssim \frac{1}{|x-y|^n}.
$$
which together with the fact that $K(y,x-y)=0$ for $|x-y|\ge d$, implies that
\eqref{K11}  holds for $N(x,y)$.
Analogously, taking derivatives in \eqref{def de K} and using that $\chi$ is a smooth function, one can show that,
$$
|\n_y N(x,y)|
\lesssim \frac{1}{|x-y|^{n+1}}
$$
which by standard arguments implies \eqref{K12}.
\end{proof}

\begin{theorem}
\label{continuidad de Tij}
The operator $T_{ij}$ is bounded on $h^p(\R^n)$ for any $\frac{n}{n+1}< p\le 1$.	
\end{theorem}
\begin{proof}
We already know that $T_{ij}$ is bounded in $L^2$. Therefore, in view of Lemma \ref{cotas del nucleo}, it only remains to show that $T_{ij}^*1\in \dot\Lambda_\a(\R^n)$ for any $0<\a<1$.	

Introducing
$$
K_1(y,z)=
\int_0^\infty\chi(y)\frac{\partial \eta_i}{\partial z_j}\left(y,\frac{z}{s}\right)\frac{ds}{s^{n+1}}
$$
and
$$
K_2(y,z)=
\int_1^\infty\chi(y)\frac{\partial \eta_i}{\partial z_j}\left(y,\frac{z}{s}\right)\frac{ds}{s^{n+1}}
$$	
we can write,
\begin{equation}
\label{desc de Tij}
\begin{aligned}
T_{ij}^*1(y)&=\lim_{\ve \rightarrow 0}\int_{|x-y|>\ve}
\chi(x) K_1(y,x-y) \, dx
-\lim_{\ve \rightarrow 0}\int_{|x-y|>\ve}
\chi(x)K_2(y,x-y) \, dx\\
&=: g_1(y) - g_2(y).
\end{aligned}
\end{equation}
For the first term observe that $K_1(y,z)$ satisfies the conditions of the classic Calder\'on-Zygmund kernels considered in \cite{CZ}. Indeed, 
if $\Sigma$ is the unit sphere,
$\int_\Sigma K_1(y,\sigma) d\sigma=0$. 
In fact,
making the change of variable $r=1/s$ we have
$$
\int_\Sigma K_1(y,\sigma) d\sigma
=\int_\Sigma\int_0^\infty \chi(y)\frac{\partial \eta_i}{\partial z_j}(y,r\sigma)
r^{n-1} dr d\sigma
=\int_{\R^n} \chi(y)\frac{\partial \eta_i}{\partial z_j}(y,z)dz =0
$$	
because $\eta_i(y, z)$ is a smooth function with compact support in the $z$ variable.
And, on the other hand, $K_1$ is homogeneous of degree $-n$ in the second variable. Indeed, making now the change of variable
$t=s/\lambda$, we obtain
$$	
K_1(y,\lambda z)
=\int_0^\infty\chi(y)\frac{\partial \eta_i}{\partial z_j}\left(y,\frac{\lambda z}{s}\right)\frac{ds}{s^{n+1}}
=\lambda^{-n}\int_0^\infty\chi(y)\frac{\partial \eta_i}{\partial z_j}\left(y,\frac{z}{t}\right)\frac{dt}{t^{n+1}}
=\lambda^{-n}K_1(y,z).
$$
Consequently, the first term on the right hand side of \eqref{desc de Tij} is a Calder\'on-Zygmund
operator acting on the smooth function $\chi(x)$, and 
so, by a known result \cite{M}, $g_1\in\dot\Lambda_\alpha(\R^n)$, for all $0<\a< 1$.
 .

On the other hand, the kernel $K_2$ is not really a singular one, indeed we have that
$$
|K_2(y,x-y)|\le
\int_1^\infty \left|\chi(y)\frac{\partial \eta_i}{\partial z_j}\left(y,\frac{x-y}{s}\right)\right|\frac{ds}{s^{n+1}}
\lesssim \int_1^\infty \frac{ds}{s^{n+1}}<\infty.
$$
Moreover $\chi(x)K_2(y,x-y)$ vanishes unless $|x-y|\le 2d$, and so 
$$
g_2(y)=\int_{|x-y|\le 2d}
\chi(x)K_2(y,x-y) \, dx.
$$
Analogously, taking derivatives with respect to $y$ in the definition of $K_2$, we can show
that $\n_y[K(y,x-y)]$ is a bounded function, and consequently,
$$
|\n g_2(y)|
\le \int_{|x-y|< d}\left|\n_y[K(y,x-y)]\right| \, dx
\lesssim 1
$$
which, since $g_2$ has compact support, implies that
$g_2\in\Lambda_\a(\R^n)$, for all $0<\a\le 1$. 
Summing up we obtain from \eqref{desc de Tij} that $T_{ij}^*1$ belongs to $\dot\Lambda_\a(\R^n)$ , for all $0<\a< 1$.

Now, given any $p\in(\frac{n}{n+1},1)$ we choose
$\a\in (0,1)$ such that $\frac{n}{n+\a}\le p\le 1$ and apply
Theorem \ref{teo de Komori Hp} to conclude the proof.
\end{proof}

We can now state the main result of this section.

\begin{theorem}
\label{teo en estrellados hp sin pesos}
Let $\O$ be a bounded domain which is star-shaped with respect to a ball $B\subset\O$. If 
$\frac{n}{n+1}<p\le 1$ and $f \in h^p_z(\O)$ satisfies $\int_\O f=0$ then, there exists $\u\in h^{1,p}_{z,0}(\O)^n$ such that
\begin{equation}
\label{solucion en estrellados en Hp}
\d\u = f \ \mbox{in} \ \O
\end{equation}
and
\begin{equation}
\label{cota en hp estrellados}
\|\u\|_{h_z^{1,p}(\O)}\le C\|f\|_{h_z^p(\O)}
\end{equation}
where $C= C(\O, p,n)$. 
\end{theorem}
\begin{proof}
By density we can assume that $f$ is smooth and so, the Bogovski\v\i{} solution of \eqref{solucion en estrellados en Hp} given by \eqref{def de u} is well defined. Moreover,  since $\sop f\subset\overline\O$, $\u$ vanishes outside $\overline\O$.

Now, it follows from \eqref{omegaij} that 
$\omega_{ij}\in C^\infty$, and therefore,
$$
\|\omega_{ij} f\|_{h_z^p(\O)}\lesssim \|f\|_{h_z^p(\O)}
$$
which combined with \eqref{derivadas de u} and Theorem \ref{continuidad de Tij} yields
$$
\left\|\frac{\partial u_i}{\partial x_j}\right\|
_{h_z^p(\O)}\lesssim \|f\|_{h_z^p(\O)}.
$$
Finally, calling $p^*=\frac{np}{n-p}$, we have
\begin{equation}
\label{cota de u en Lp estrella}
\|\u\|_{h_z^p(\O)}\lesssim C\|\u\|_{L^{p^*}(\O)}
\lesssim\|\n\u\|_{h_z^p(\O)}
\end{equation}
where we have used, for the first inequality  that $L^{p^*}(\O)\subset h^p(\O)$, and for the second one
a theorem from \cite{KS}. Although the proof of this embedding given in that paper is for $H^p(\R^n)$ it is not difficult to see that the same arguments can be applied to $h^p(\O)$.
\end{proof}

As a consequence of the previous theorem we have the following result which can be proved by known arguments.

\begin{corollary}
\label{dominios lipschitz}
Let $\O$ be a bounded Lipschitz domain. If 
$\frac{n}{n+1}<p\le 1$ and $f\in h^p_z(\O)$ then, there exists $\u\in h^{1,p}_{z,0}(\O)^n$ such that
$$
\d\u = f \ \mbox{in} \ \O
$$
and
$$
\|\u\|_{h_z^{1,p}(\O)}\le C\|f\|_{h_z^p(\O)}
$$
where $C= C(\O, p,n)$. 
\end{corollary}

\begin{proof}
$\O$ can be decomposed into a finite union of star-shaped domains $\O_1, \O_2, \cdots,\O_m$, and any $f\in h^p_z(\O)$ can be written as
$f=\sum_{i=1}^m f_i$ , $\int f_i=0$, $f_i\in h^p_z(\O_i)$, $\|f_i\|_{h^{p}_z(\O)}\lesssim C\|f\|_{h_z^p(\O)}$, where for these inequalities we have to use \eqref{eta por f esta en hp}. The existence of this decomposition can be proved by the same arguments used in \cite[Pag. 132, Lemma 3.4]{Ga}. Now, if $\u_i$ is the Bogovski\v i's solution of $\d\u_i=f_i$ in $\O_i$ then $\u=\sum_{i=1}^m\u_i$ is the desired solution in $\O$.
\end{proof}

\section{Generalization to weighted Hardy spaces}

The results of the previous section can be extended without too much difficult to the weighted case.

Recall that given a non-negative locally integrable $w$, the weighted spaces $H_w^p(\R^n)$ and $h_w^p(\R^n)$ are defined as (see for example \cite{StroTor}),
$$
H_w^p(\R^n)=\{f\in\mathcal S'\colon M_{\varphi}f\in L_w^p(\R^n)\},
$$
and 
$$
h_w^p(\R^n)=\{f\in\mathcal S'\colon m_{\varphi}f\in L_w^p(\R^n)\}.
$$
With obvious generalization for the spaces on domains and their notation, we have the following result. 

\begin{theorem}
Let $\O$ be a bounded domain which is star-shaped with respect to a ball $B\subset\O$, and assume that $w\in A_1$. If 
$\frac{n}{n+1}<p\le 1$ and $f\in h^p_{w,z}(\O)$
satisfies $\int_\O f=0$ then, there exists 
$\u\in h^{1,p}_{z,w,0}(\O)^n$ such that
$$
\d\u = f \ \mbox{in} \ \O
$$
and
$$
\|\u\|_{h_{z,w}^{1,p}(\O)}\le C\|f\|_{h_{z,w}^p(\O)}
$$
where $C= C(\O, p,n,w)$. 
\end{theorem}
\begin{proof} 
The same arguments used in Theorem \ref{teo en estrellados hp sin pesos} can be applied to estimate
the derivatives of $\u$. 
The only important difference is that we have to use a weighted version of Theorem \ref{teo de Komori Hp}
which follows easily from 
\cite{K3} by the same argument explained in 
that theorem. Finally, a weighted version of the embedding \eqref{cota de u en Lp estrella}
follows from 
\cite[Theorem 8.7]{H} applied to $\R^n$ with measure $\mu=wdx$.
\end{proof}

\section{Solvability of the divergence in Lipschitz spaces}
\label{seccion lipschitz}

As in the previous section we consider first the case of star-shaped domains and the solution given in \eqref{def de u}.

To prove the continuity of $T_{ij}$ on Lipschitz spaces
we will make use of the following known result.

\begin{theorem} 
\label{teo de Komori lipschitz}
Let 
$$
Tf(x)=\lim_{\ve \rightarrow 0}\int_{|x-y|>\ve} N(x,y)f(y)\,dy
$$
be a singular integral operator which is bounded in $L^2$. 
Assume that $T$ satisfies the following conditions:
	
\begin{equation}
\label{K11 Lipschitz}
|N(x,y)|\lesssim \min\left\{\frac{1}{|x-y|^n},\frac{1}{|x-y|^{n+1}}\right\}   
\end{equation}
	
\begin{equation}
\label{K12 Lipschitz}
2|x-\bar x| < |y - x| \implies |N(x, y) - N(\bar x,y)| \lesssim \frac{|x-\bar x|}{|x-y|^{n+1}}.
\end{equation}
If, for all $0<\a<1$, $T1\in\dot\Lambda_\a(\R^n)$, then 
$T$ is bounded on $\Lambda_\a(\R^n)$, for all $0<\a<1$, as well as on $bmo(\R^n)$.
\end{theorem} 
\begin{proof}
Assuming $|N(x,y)|\lesssim |x-y|^{-n}$ instead of \eqref{K11 Lipschitz} Komori considered singular integral operators in more general spaces. He proved a result that implies in particular that $T$ is continuous from $\Lambda_\a(\R^n)$ to
$\dot\Lambda_\a(\R^n)$, for $0<\a<1$ (see \cite[Corollary 1]{K2}). His arguments could be modified to show that, assuming \eqref{K11 Lipschitz}, $T$ is continuous on $\Lambda_\a(\R^n)$.

However, the weaker result stated here follows immediately by duality. Indeed, our hypotheses are exactly
those in Theorem \ref{teo de Komori Hp} but interchanging variables $x$ and $y$. Therefore, it follows from that theorem that $T^*$ is continuous on $h^p(\R^n)$ for all $\frac{n}{n+1}<p\le 1$. 

Thus, the proof concludes recalling that
$(h^1(\R^n))^*=bmo(\R^n)$, $(h^p(\R^n))^*=\Lambda_{n(1/p-1)}(\R^n)$, and given $\a\in(0,1)$, applying the continuity of $T^*$ in $h^{\frac{n}{n+\a}}(\R^n)$.
\end{proof}

To simplify notation, the continuity of $T_{ij}$ will be reduced below to the continuity of
operators of the kind analyzed in the following lemma.

\begin{lemma}
\label{lema H}
Let 
$$
Tf(x)=\lim_{\ve \rightarrow 0}\int_{|x-y|>\ve} H(x,y)f(y)\,dy
$$
be a singular integral operator
with kernel
\begin{equation}
\label{H}
H(x,y)=\int_0^1 \chi(x) \chi(y)
\psi\left( y+\frac{x-y}{s}\right) \frac{ds}{s^{n+1}},
\end{equation}
where $\psi\in C^{\infty}_0(B)$ has zero integral.

Then, 
$T$ is bounded on $\Lambda_\a(\R^n)$, for all $0<\a<1$, as well as on $bmo(\R^n)$.
\end{lemma}

\begin{proof}
We decompose $H(x,y)$ as,

\begin{align}
H(x,y)&=\int_0^{1/2} \chi(x) \chi(y)
\psi\left( y+\frac{x-y}{s}\right) \frac{ds}{s^{n+1}}
+\int_{1/2}^1 \chi(x) \chi(y)
\psi\left( y+\frac{x-y}{s}\right) \frac{ds}{s^{n+1}}\\
&=: H_1(x,y)+H_2(x,y)
\end{align}

and consider the corresponding operators $T_1$ and $T_2$.

The operator $T_2$ is not a singular one and so it is easy to analyze. Indeed, $H_2(x,y)$ is a smooth  function with all its derivatives bounded, and therefore,
$T_2$ is bounded on $\Lambda_\a(\R^n)$, for  $0<\a<1$, and on $bmo(\R^n)$.

For $T_1$ writing 
$$
\psi\left( y+\frac{x-y}{s}\right) 
=\psi\left( x+\frac{x-y}{s}-(x-y)\right) 
$$
we have,
$$
T_1f(x)=\lim_{\ve \rightarrow 0}\int_{|x-y|>\ve}\chi(y) K(x,x-y) f(y) \,dy
$$
where
$$
K(x,z)=\int_0^{1/2} 
\chi(x)\psi\left( x+\frac{z}{s}-z\right) \frac{ds}{s^{n+1}},
$$
and making the change of variable $t=s/(1-s)$ we obtain,
$$
K(x,z)=
\int_0^1 \chi(x)
\psi\left( x+\frac{z}{t}\right) (1+t)^{n-1}\frac{dt}{t^{n+1}}.
$$
By the same arguments used in Lemma \ref{cotas del nucleo} and Theorem \ref{continuidad de Tij} we can see
that the kernel $N(x,y): = K(x,x-y)$ satisfies the bounds
\eqref{K11 Lipschitz} and \eqref{K12 Lipschitz}

To apply Theorem \ref{teo de Komori lipschitz} to  $T_1$
it only remains to prove that $T_1$ is bounded on $L^2$ and that  $T_11\in\dot\Lambda_\a(\R^n)$, for all $0<\a<1$.

Expanding $(1+t)^{n-1}$ by the Newton binomial formula we
obtain,
$$
K(x,z)=K_0(x,z)+\cdots+K_{n-1}(x,z)
$$
where 
$$
K_0(x,z)=
\int_0^1 \chi(x)
\psi\left( x+\frac{z}{t}\right) \frac{dt}{t^{n+1}}
$$ 
and for $j=1,\cdots, n-1$,
\begin{equation}
\label{nucleo de Sj}
K_j(x,z)=\binom{n-1}{j}
\int_0^1 \chi(x)
\psi\left( x+\frac{z}{t}\right) \frac{dt}{t^{n+1-j}}.
\end{equation}

Consider the associated decomposition
$$
T_1=S_0+ \cdots + S_{n-1}
$$
with
$$
S_jf(x)=\lim_{\ve \rightarrow 0}\int_{|x-y|>\ve} \chi(y)K_j(x,x-y) f(y) \,dy.
$$
Now, to prove that $S_01\in\dot\Lambda_\a(\R^n)$  we decompose $K_0$ as
\begin{equation}
\label{S0}
K_0(x,z)=
\int_0^\infty \chi(x)
\psi\left( x+\frac{z}{t}\right) \frac{dt}{t^{n+1}}
-\int_1^\infty \chi(x)
\psi\left( x+\frac{z}{t}\right) \frac{dt}{t^{n+1}}
\end{equation}
and consider the two associated operators. The first one satisfies the conditions of the classic Calder\'on-Zygmund operators analyzed in \cite{CZ},
while the second one has a smooth kernel,
and consequently, 
$S_0$ is bounded on $L^2$ and $S_01\in\dot\Lambda_\a(\R^n)$.
We omit details because they are exactly as those used in
Theorem \ref{continuidad de Tij} for $T^*_{ij}1$.

It remains to analyze $S_j$ for $j=1,\cdots,n-1$. These are not really singular integral operators, indeed, using that the integrand in \eqref{nucleo de Sj} vanishes unless $|z|<td$, we have
$$
|K_j(x,x-y)|=\binom{n-1}{j}
\int_{|x-y|/d}^1 \chi(x)
\left|\psi\left( x+\frac{x-y}{t}\right)\right| \frac{dt}{t^{n+1-j}}
\lesssim \frac{1}{|x-y|^{n-j}}.
$$
which implies immediately, by the H\"older inequality
that $S_j$ is bounded on $L^2$.

On the other hand, given $x$ and $\bar x$ we write
$$
\begin{aligned}
S_j1(x)&-S_j1(\bar x)
=\int_{|x-y|<|\bar x-y|} \chi(y)[K_j(x,x-y)-K_j(\bar x,\bar x-y)] \, dy \\
&+\int_{|\bar x-y|<|x-y|} \chi(y)[K_j(x,x-y)-K_j(\bar x,\bar x-y)] \, dy =: I + II.
\end{aligned}
$$
For $|x-y|<|\bar x-y|$ we have,
for any $0<\a<1$,
$$
\begin{aligned}
|K_j(x,x-y)-K_j(\bar x,\bar x-y)|
&\lesssim
\int_{|x-y|/d}^1 \left|\chi(x)
\psi\left( x+\frac{x-y}{t}\right)
-\chi(\bar x)
\psi\left(\bar x+\frac{\bar x-y}{t}\right)
\right| \frac{dt}{t^{n+1-j}} \\
&\lesssim |x-\bar x|^\a\int_{|x-y|/d}^1 
\frac{dt}{t^{n+1-j+\a}}
\lesssim \frac{|x-\bar x|^\a}{|x-y|^{n-j+\a}}
\end{aligned}
$$
where we have used that $\chi$ and $\psi$ are smooth functions with compact support.

Consequently, since $j\ge 1$ and $\a<1$, we obtain
$$
|I|\lesssim \int_{|x-y|<d} \frac{|x-\bar x|^\a}{|x-y|^{n-j+\a}} \, dy
\lesssim |x-\bar x|^\a.
$$
Analogously we can estimate $II$, and therefore,
$S_j1\in\dot\Lambda_\a(\R^n)$, for all $0<\a<1$, as we wanted to prove.
\end{proof}

\begin{theorem}
\label{teo caso Lipschitz}
Let $\O$ be a bounded domain which is star-shaped with respect to a ball $B\subset\O$. 
If $f\in\Lambda_\a(\R^n)$, for $0<\a<1$ (resp. $f\in bmo(\R^n)$) is such that
$\sop f\subset\overline\O$
and $\int_\O f=0$ then, there exists $\u\in\Lambda_\a(\R^n)^n$ (resp. $\u\in bmo(\R^n)^n$) such that 
$\n\u\in\Lambda_\a(\R^n)^{n\times n}$ (resp. $\n\u\in bmo(\R^n)^{n\times n}$), 
$\sop\u\subset\overline\O$, and
$$
\d\u = f \ \mbox{in} \ \O
$$
and
\begin{equation}
\label{cota de las derivadas en Lipschitz alpha}
\left\|\frac{\partial u_i}{\partial x_j}\right\|
_{\Lambda_\a(\R^n)}
\le C\|f\|_{\Lambda_\a(\R^n)} 
\end{equation}
where $C= C(\O,\a,n)$, or 
\begin{equation}
\label{cota de las derivadas en bmo}
\left\|\frac{\partial u_i}{\partial x_j}\right\|
_{bmo(\R^n)}
\le C\|f\|_{bmo(\R^n)}
\end{equation}
where $C= C(\O,n)$.
\end{theorem}
\begin{proof}
It is enough to prove \eqref{cota de las derivadas en Lipschitz alpha} and \eqref{cota de las derivadas en bmo}.
Let $\u$ be the solution defined in \eqref{def de u}. Writing
$$
\begin{aligned}
\frac{\partial G}{\partial x_j}(x,y)
&=\int_0^1\frac{\partial}{\partial x_j}
\left[\frac{x-y}{s}\omega\left( y+\frac{x-y}{s}\right) 
\right] \frac{ds}{s^n}\\
& = \int_0^1\frac{\partial}{\partial x_j}
\left[\left(y+\frac{x-y}{s}\right)\omega\left( y+\frac{x-y}{s}\right)\right] \frac{ds}{s^n}
-\int_0^1 y\frac{\partial}{\partial x_j}
\left[\omega\left( y+\frac{x-y}{s}\right)\right] \frac{ds}{s^n},
\end{aligned}
$$
we have
\begin{equation}
\label{descomposicion para simplificar}
T_{ij}f(x)=T_1f(x) - T_2 (yf)(x)
\end{equation}
where $T_1$ and $T_2$ are singular integral operators
with kernel of the form
$$
H(x,y)=\int_0^1 \chi(x)\chi(y)\frac{\partial}{\partial x_j}	\left[\varphi\left( y+\frac{x-y}{s}\right)\right] \frac{ds}{s^n}
$$
with $\varphi=x\omega$ for $T_1$ and $\varphi=\omega$ for $T_2$. Thus, both operators are of the form
given in \eqref{H} with $\psi=\partial\vphi/\partial x_j$. Thus, \eqref{cota de las derivadas en Lipschitz alpha} and \eqref{cota de las derivadas en bmo} follow from the representation \eqref{derivadas de u}, recalling that $\omega_{ij}$ is a smooth function and applying Lemma \ref{lema H}.
\end{proof}

As in the case of Hardy spaces we can extend the result to Lipschitz domains.

\begin{corollary}
\label{dominios lipschitz para Lipschitz alpha}
Let $\O$ be a bounded Lipschitz domain. 
If $f\in\Lambda_\a(\R^n)$, for $0<\a<1$ (resp. $f\in bmo(\R^n)$) is such that
$\sop f\subset\overline\O$
and $\int_\O f=0$ then, there exists $\u\in\Lambda_\a(\R^n)^n$ (resp. $\u\in bmo(\R^n)^n$) such that 
$\n\u\in\Lambda_\a(\R^n)^{n\times n}$ (resp. $\n\u\in bmo(\R^n)^{n\times n}$), 
$\sop\u\subset\overline\O$, and
$$
\d\u = f \ \mbox{in} \ \O
$$
and
$$
\left\|\frac{\partial u_i}{\partial x_j}\right\|
_{\Lambda_\a(\R^n)}
\le C\|f\|_{\Lambda_\a(\R^n)} 
$$
where $C= C(\O,\a,n)$, or 
$$
\left\|\frac{\partial u_i}{\partial x_j}\right\|
_{bmo(\R^n)}
\le C\|f\|_{bmo(\R^n)}
$$
where $C= C(\O,n)$.
\end{corollary}
\begin{proof} The result follows by the arguments given in \cite[Pag. 132, Lemma 3.4]{Ga} that we have already mentioned in Corollary \ref{dominios lipschitz}.
\end{proof}

\begin{remark} The condition 
$\sop f\subset\overline\O$ is essential to have the estimate \eqref{cota de las derivadas en Lipschitz alpha} for $\u$ given by \eqref{def de u}. Indeed, that solution satisfies $\d\u=f$ in $\R^n$ and vanishes outside $\O$, therefore, $\n\u\in\Lambda_\a(\R^n)^{n\times n}$ implies that
$f$ vanishes outside $\O$. However, the problem \eqref{div} has infinitely many solutions, and therefore, it is an interesting problem whether there exists a solution having first derivatives
belonging to $\Lambda_\alpha(\O)$ for arbitrary $f\in \Lambda_\alpha(\O)$ satisfying the zero integral condition. Some result in this direction was obtained by Berselli and Longo \cite{BL} who
proved, under appropriate regularity assumptions on the domain, that there exists a $C^1(\bar\O)$ solution if $f$ satisfies a Dini condition (in particular, if $f\in\Lambda_\alpha(\O)$).
\end{remark}

\section{The Korn inequality}

It is well known that the different classic versions of the Korn inequality for vector fields $\u\in W^{1,p}(\O)^n$, for $1<p<\infty$, can be derived from
Theorem \ref{teorema en Lp} (see for example \cite{AD}). In this section we extend the Korn inequality to the case $\frac{n}{n+1}<p\le 1$ using the continuity of the singular integral operators analyzed in the previous sections.

With this goal we first recall that the strain tensor $\ve(\u)$ of a given vector field $\u$ is defined as the symmetric part of its differential matrix, namely,
$$
\ve_{ij}(\u)=\frac12 \left(\frac{\partial u_i}{\partial x_j} + \frac{\partial u_j}{\partial x_i}\right)
$$
To simplify notation we drop the $\u$ and write simply $\ve_{ij}$ for $\ve_{ij}(\u)$ and $\ve_{ij,k}$ for
$\frac{\partial{\ve_{ij}}}{\partial{x_k}}$. 

\begin{theorem} Let $\O$ be a bounded Lipschitz domain. For $\frac{n}{n+1}<p\le 1$, if $\u\in h^{1,p}_r(\O)$ then,
\begin{equation}
\label{korn en hp}
\|\nabla\u\|_{h^p_r(\O)}
\le C \left\{\ve(\u)\|_{h^p_r(\O)} + \|\u\|_{h^p_r(\O)}\right\}
\end{equation}
where $C(\O,n,p)$.	
\end{theorem}

\begin{proof}
We assume that $\O$ is star-shaped with respect to a ball $B\subset\O$. The general case will follow
using again that a bounded Lipschitz domain can be written as a finite union of star-shaped domains.

Given $\vphi\in C^\infty(\R^n)$ one can show by simple computations (see \cite{AD} for details) that, for  $y\in\O$,
\begin{equation}
\label{representacion de fi}
\vphi(y)-\vphi_\omega
=-\int_\O \G(x,y)\cdot \n\vphi(x) \, dx
\end{equation}
where $\vphi_\omega=\int \vphi\omega$ and $\G(x,y)$ is the kernel defined in \eqref{definicion de G}.

By density we can assume that $\u$ is smooth, and then, applying \eqref{representacion de fi} with $\vphi=\frac{\partial u_i}{\partial x_j}$, we have
$$
\frac{\partial u_i}{\partial x_j}(y)
-\left(\frac{\partial u_i}{\partial x_j}\right)_\omega
=-\sum_{k=1}^n\int_\O G_k(x,y)
\frac{\partial^2 u_i}{\partial x_k\partial x_j}(x)\, dx.
$$	
which, using the elementary identity
$$
\frac{\partial^2 u_i}{\partial x_k\partial x_j}
=\ve_{ik,j}
+\ve_{ij,k}
- \ve_{jk,i},
$$
becomes
$$
\frac{\partial u_i}{\partial x_j}(y)
-\left(\frac{\partial u_i}{\partial x_j}\right)_\omega
=-\sum_{k=1}^n\int_\O G_k(x,y)
(\ve_{ik,j}(x)
+\ve_{ij,k}(x)
- \ve_{jk,i}(x))
\, dx.
$$
for any $y\in\O$.

Now, recalling that $G_k(x,y)=0$ for all $y\in\O$ and $x\in\partial\O$ we can integrate by parts in the usual way for singular integrals (see details in \cite[Lemma 2.3]{AD}) to obtain, for $y\in\O$ and $x\in\O$,
\begin{equation}
\label{final de Korn}
\begin{aligned}
\frac{\partial u_i}{\partial x_j}(y)
-\left(\frac{\partial u_i}{\partial x_j}\right)_\omega
&=\sum_{k=1}^n
\left(T^*_{kj}\ve_{ik}(y)+T^*_{kk}\ve_{ij}(y)
-T^*_{ki}\ve_{jk}(y)\right)\\
&+\sum_{k=1}^n
\left(\omega_{kj}(y)\ve_{ik}(y)
+\omega_{kk}(y)\ve_{ij}(y)
-\omega_{ki}(y)\ve_{jk}(y)\right)
\end{aligned}
\end{equation}
where, for $i,j=1,\cdots,n$, $\omega_{ij}$ and $T_{ij}$
are as in \eqref{omegaij} and \eqref{operador general}.
Note that we are considering now $\ve_{ij}$ extended to $h^p_r(\R^n)$, which can be done because $\u\in h^{1,p}_r(\O)$.

Now, to estimate the first line on the right hand side of \eqref{final de Korn} we use Theorem \ref{teo de Komori Hp} applied now to $T^*_{ij}$. Indeed, all the hypotheses of that theorem have been verified in Section \ref{seccion lipschitz}. 

On the other hand, for the second line of \eqref{final de Korn} we recall that $\omega_{ij}$ are smooth functions and 
use \eqref{eta por f esta en hp}.

Finally,
$$
\left|\left(\frac{\partial u_i}{\partial x_j}\right)_\omega\right|
=\left|\int \frac{\partial u_i}{\partial x_j}(x)\,\omega(x)dx\right|
\le \left|\int u_i(x)\frac{\partial \omega}{\partial x_j}(x)
dx\right|
\lesssim C_{\omega,p}
\|u_i\|_{h^p_r(\O)}.
$$
where
$C_{\omega,p}=\left\|
\frac{\partial \omega}{\partial x_j}\right\|_{\Lambda_{n(1/p-1)}(\R^n)}$,
for $\frac{n}{n+1}<p<1$, while
$C_{\omega,1}=\left\|
\frac{\partial \omega}{\partial x_j}\right\|_{bmo(\R^n)}$.
\end{proof}

\end{document}